\documentclass[reqno,12pt]{amsart}
\usepackage{amssymb}

\theoremstyle{plain}    
\newtheorem{thm}{Theorem}[section]
\numberwithin{figure}{section} 
\theoremstyle{plain}    
\newtheorem*{thm*}{Theorem} 
\theoremstyle{plain}    
\newtheorem{cor}[thm]{Corollary} 
\theoremstyle{plain}    
\newtheorem{lem}[thm]{Lemma} 
\theoremstyle{plain}    
\newtheorem{prop}[thm]{Proposition} 
\theoremstyle{definition}

\theoremstyle{remark}

\theoremstyle{remark}

\theoremstyle{remark}    

\theoremstyle{remark}    

\theoremstyle{definition}  
\newtheorem{example}[thm]{Example}
\theoremstyle{remark}
  \newtheorem*{acknowledgement*}{Acknowledgement} 

\theoremstyle{plain}    

\theoremstyle{plain}    
\theoremstyle{plain}    
\theoremstyle{plain}    
\theoremstyle{definition}

\theoremstyle{remark}

\theoremstyle{remark}    

\theoremstyle{remark}    

\theoremstyle{plain}    

\newcount\theTime
\newcount\theHour
\newcount\theMinute
\newcount\theMinuteTens
\newcount\theScratch
\theTime=\number\time
\theHour=\theTime
\divide\theHour by 60
\theScratch=\theHour
\multiply\theScratch by 60
\theMinute=\theTime
\advance\theMinute by -\theScratch
\theMinuteTens=\theMinute
\divide\theMinuteTens by 10
\theScratch=\theMinuteTens
\multiply\theScratch by 10
\advance\theMinute by -\theScratch

\def\today{{\number\day\space
 \ifcase\month\or
  January\or February\or March\or April\or May\or June\or
  July\or August\or September\or October\or November\or December\fi
 \space\number\year}}

\newcommand\Afr{{\mathfrak A}}

\newcommand\alphat{{\tilde\alpha}}

\newcommand\betat{{\tilde\beta}}

\newcommand\Cpx{{\mathbf C}}

\newcommand\Dfr{{\mathfrak D}}

\newcommand\dist{\operatorname{dist}}

\newcommand\eps{\epsilon}

\newcommand\Et{{\widetilde E}}

\newcommand\Fb{{\mathbf F}}

\newcommand\HEu{\mathcal{H}}

\newcommand\Ints{{\mathbf Z}}

\newcommand\Lambdao{{\Lambda\oup}}

\newcommand\lspan{\mathrm{span}\,}

\newcommand\Mcal{{\mathcal{M}}}

\newcommand\Nats{{\mathbf N}}

\newcommand\oup{^{\mathrm o}}

\newcommand\pit{{\tilde\pi}}

\newcommand\rank{\mathrm{rank}\,}

\newcommand\restrict{{\upharpoonright}}

\newcommand\tr{{\operatorname{tr}}}

\newcommand\dif{\mbox{\it d}}

\newcommand\gammat{{\tilde\gamma}}

\newcommand\Reals{{\mathbf R}}

\newcommand\vt{{\tilde v}}

\begin{document}

\pagestyle{myheadings}

\title{Popa algebras in free group factors}

\date{19 June, 2002}

\author{Nathanial P. Brown}
\address{Department of Mathematics, Michigan State University, 
East Lansing, MI 48824}
\email{brown1np@cmich.edu}

\author{Kenneth J. Dykema} 
\address{Department of Mathematics, Texas A\&M University, 
College Station, TX 77843} 
\email{kjd@tamu.edu}

\thanks{N.B.\ is an NSF Postdoctoral Fellow. K.D.\ is supported in 
part by NSF grant DMS--0070558.}

\begin{abstract}
For each $1<s<\infty$, a Popa algebra $A_s$ is constructed
that embeds as a weakly dense C$^*$--subalgebra of the interpolated free group factor
$L(\Fb_s)$.
Certain approximation properties for $A_s$ are shown.
It follows that $L(\Fb_s)$ has the weak expectation property of Lance
with respect to $A_s$.
In the course of the demonstration, it is proved that under certain conditions,
full amalgamated free products of matrix algebras are residually finite dimensional.
\end{abstract}
\maketitle

\markboth{Popa Algebras}{Popa Algebras}

\section{Introduction}

The interpolated free group factors $L(\Fb_s)$, $1<s\le\infty$,
(see~\cite{Ra94} and~\cite{D:interp}) are a family of von Neumann
algebra II$_1$--factors including the usual free group factors (when
$s\in\Nats\cup\{\infty\}$).  In this paper, we construct for every
$1<s<\infty$ a weakly dense, unital C$^*$--subalgebra $A_s\subseteq
L(\Fb_s)$ which has the following properties:
\renewcommand{\labelenumi}{(\Roman{enumi})}
\begin{enumerate}
\item $A_s$ is a finitely generated Popa algebra,
\item $L(\Fb_s)$ has the weak expectation property (WEP) of Lance~\cite{lance:nuclear}
      relative to $A_s$.
\end{enumerate}

By definition, a Popa algebra is a unital, simple C$^{*}$-algebra $A$ with 
the property that for every finite set $\mathfrak{F} \subset A$ and 
$\epsilon > 0$ there exists a finite dimensional subalgebra $0 \neq B 
\subset A$ (the inclusion not necessarily unital) with unit $e$ such 
that $\| [x,e] \| < \epsilon$ for all $x \in \mathfrak{F}$ and for 
each $x \in \mathfrak{F}$ there exists $b \in B$ such that
$\| b -  exe \| < \epsilon$.

It was shown in~\cite{brown:AFDtraces} that all McDuff II$_1$--factors
have weakly dense Popa algebras, but with Property~(I) above, our
construction yields the first non--McDuff example, which thus answers
Question~11.1 in~\cite{brown:AFDtraces}.  Since $A_s$ is a Popa
algebra, it is quasidiagonal (see~\cite{popa:simpleQD}, \cite{brown:QDsurvey}
and~\cite{dvv:QDsurvey}) and thus quite distinct from previously known
weakly dense C$^*$--subalgebras of $L(\Fb_s)$ that arise as reduced
free products of C$^*$--algebras.  Indeed, since $C^*_r ({\bf
F}_2)$ is not quasidiagonal and quasidiagonality passes to
subalgebras, it follows that no C$^*$-algebra containing a copy of
$C^*_r ({\bf F}_2)$ (e.g.\ most reduced free products) can be
quasidiagonal.
(Moreover, other ways in which $A_s$ differs from reduced free product C$^*$--algebras
are discussed at the end of~\S6.)
Thus we regard the algebras $A_{s}$ as rather exotic
weakly dense subalgebras of $L(\Fb_s)$ and we believe they provide
interesting examples to test the invariance questions around
Voiculescu's free entropy (see~\cite{V94} and~\cite{dvv:partIII}).
Indeed, if the free entropy dimension of a generating set for $A_s$ were
computed and found to be different from $s$, this would show that free entropy
dimension of a generating set is not an invariant for von Neumann algebras.
On the other hand, if the free entropy dimension were found to be $s$,
then because $A_s$ is so different from previously considered
examples, this would be nontrivial evidence for invariance.

The meaning of Property~(II) above is that for any normal, unital
representation $L(\Fb_s)\subset B(H)$ with $H$ a separable Hilbert
space, there is an idempotent, completely positive map $\Phi:B(H)\to
L(\Fb_s)$ such that $\Phi(a)=a$ for all $a\in A_s$.  Since $A_s$ lies
in the multiplicative domain of $\Phi$ (cf.\ \cite{paulsen:cbmaps}) it
follows that we also have the bi-module property $\Phi(axb)=a\Phi(x)b$
for all $a,b\in A_s$ and $x\in B(H)$.  It was shown
in~\cite{brown:AFDtraces} that a McDuff factor has the WEP relative to
a weakly dense C$^*$--subalgebra if and only if it approximately
embeds into the hyperfinite II$_1$--factor.  However, with Property~(II),
we give the first example of non-McDuff
II$_1$--factors with a weak expectation.

The algebras $A_{s}$ are constructed as (corners of) inductive limits
of certain full amalgamated free products of finite dimensional
C$^*$-algebras.  This construction is an adaptation of the techniques
in~\cite{brown:AFDtraces} which apply some basic ideas from Elliott's
classification program.  We will show that each $A_s$ has a trace
whose GNS representation generates a von Neumann algebra that is
isomorphic to a certain amalgamated free product of hyperfinite
II$_1$--factors.  Using some standard techniques from free probability
theory (cf.~\cite{D:amalg} or~\cite{D:subf}), we will prove that these
amalgamated free products are, in turn, isomorphic to interpolated
free group factors.  The first named author thanks Dima Shlyakhtenko
for first suggesting that these amalgamated free products would be
isomorphic to interpolated free group factors.

This paper is organized as follows.  In~\S2,
we show that certain full amalgamated free products of
matrix algebras are residually finite dimensional.  This result will
be needed to ensure that the inductive limits we consider give simple
C$^{*}$-algebras.  In~\S3 we prove the result described above about
certain amalgamated free products of von Neumann algebras being
interpolated free group factors.  In~\S4, we construct the Popa
algebras $A_s$ and traces on them, proving that their GNS
representations yield interpolated free group factors.  In~\S5 we
present a technical result which is needed to prove Property~(II)
above.  Namely we show that the canonical free product traces on
certain full amalgamated free products are weakly approximately finite
dimensional in the sense of \cite[Definition 3.1]{brown:AFDtraces}.
The proof involves a lifting result for representations of
finite dimensional C$^*$--algebras in ultraproducts of
II$_1$--factors.  Finally, in~\S6 we prove some approximation
properties of the weakly dense subalgebras $A_s$; these
imply that $L(\Fb_s)$ has the WEP relative to $A_s$.

We will use standard notation and terminology in operator algebra 
theory throughout the paper.  In particular, the symbols $\otimes$ and   
$\bar{\otimes}$ will always denote the spatial C$^{*}$ and W$^{*}$ 
tensor products, respectively. 

\section{Amalgamated Free Products of Finite Dimensional C$^*$-algebras} 

Recall that a C$^*$--algebra is \emph{residually finite dimensional} (r.f.d.) if it has a separating family of
finite dimensional representations.
The main result of this section is that the full C$^*$--algebra amalgamated free product of matrix algebras is
a residually finite dimensional C$^*$--algebra, 
provided that the normalized traces on the matrix algebras restrict to the same trace
on the subalgebra over which one amalgamates.
(That this trace condition cannot be dropped is easy to see --- cf.\ Example~\ref{ex:notrfd}.)
This result looks back to Exel and Loring's result~\cite{EL}
that the full C$^*$--algebra free product of residually finite dimensional C$^*$--algebras with amalgamation
over the scalars is residually finite dimensional.

We will use the following lemma, which is well known. 

\begin{lem}
\label{lem:pAp}
Let $A$ be a unital 
C$^*$--algebra and suppose there is a projection $p\in A$ and there are partial isometries
$v_1,\ldots,v_n\in A$ such that $v_i^*v_i\le p$ and $\sum_{i=1}^nv_iv_i^*=1-p$.
Let $v_0=p$.
\renewcommand{\labelenumi}{(\roman{enumi})}
\begin{enumerate}
\item If $S\subseteq A$ generates $A$ as a C$^*$--algebra then
\[
\bigcup_{i,j=0}^nv_i^*Sv_j
\]
generates $pAp$ as a C$^*$--algebra.
\item A is r.f.d.\ if and only if $pAp$ is r.f.d.
\end{enumerate}
\end{lem}
\begin{proof}
For~(i), take any word $x=x_1x_2\cdots x_k$
in $S\cup S^*$ and interpose $1=v_0v_0^*+v_1^*v_1+\cdots+v_n^*v_n$
between all letters $x_jx_{j+1}$ in this word.
Now compress by $p$.

For~(ii), it is clear from the definition that $A$ r.f.d.\ implies $pAp$ r.f.d., for any projection $p\in A$.
Now suppose $pAp$ is r.f.d.
Then clearly $M_{n+1}(pAp)$ is r.f.d.
The map $A\to M_{n+1}(pAp)$ given by $a\mapsto(v_i^*av_j)_{i,j=0}^n$ is a $*$--monomorphism identifying
$A$ with a corner of $M_{n+1}(pAp)$.
\end{proof}

The following lemma is certainly well known; however, for completeness we give a proof.

\begin{lem}
\label{lem:fpcorner}
Let $A$ and $B$ be unital C$^*$--algebras having a C$^*$--algebra $D$ embedded as a unital C$^*$--subalgebra
of each of them and let
\[
\Afr=A*_DB
\]
be the full C$^*$--algebra free product of $A$ and $B$ with amalgamation over $D$.
Suppose there is a projection $p\in D$ and there are partial isometries $v_1,\ldots,v_n\in D$ such that
$v_i^*v_i\le p$ and $\sum_{i=1}^nv_iv_i^*=1-p$.
Then
\[
p\Afr p\cong (pAp)*_{pDp}(pBp)\;.
\]
More precisely, if $\iota_A:A\to\Afr$ and $\iota_B:B\to\Afr$ are the inclusions arising in the free product construction,
then there is an isomorphism
\begin{equation}
\label{eq:piso}
p\Afr p\overset{\sim}{\longrightarrow}(pAp)*_{pDp}(pBp)
\end{equation}
intertwining the inclusion $\iota_A\restrict_{pAp}:pAp\to p\Afr p$, respectively
$\iota_B\restrict_{pBp}:pBp\to p\Afr p$, with the inclusion of $pAp$, respectively $pBp$,
into the RHS of~\eqref{eq:piso} arising from the free product construction.
\end{lem}
\begin{proof}
From part~(i) of Lemma~\ref{lem:pAp}, $p\Afr p$ is generated by $pAp\cup pBp$.
It remains to show that if $\HEu$ is a Hilbert space and
\begin{align*}
\pi:pAp&\to B(\HEu) \\
\sigma:pBp&\to B(\HEu)
\end{align*}
are $*$--homomorphisms such that $\pi\restrict_{pDp}=\sigma\restrict_{pDp}$,
then there is a $*$--homomorphism $\psi:p\Afr p\to B(\HEu)$ such that
\begin{align*}
\psi\circ\iota_A\restrict_{pAp}&=\pi \\
\psi\circ\iota_B\restrict_{pBp}&=\sigma\;.
\end{align*}
Let $\HEu_i$ be a copy of $\pi(v_i^*v_i)\HEu$;
also let $\HEu_0=\HEu$ and $v_0=p$.
Let
\[
\HEu'=\bigoplus_{i=0}^n\HEu_i\;.
\]
We will correspondingly think of an operator $T\in B(\HEu')$ as being an $(n+1)\times(n+1)$ matrix
whose $(i,j)$th entry $T_{ij}$,
($0\le i,j\le n$), is the part of $T$ sending direct summand $\HEu_j$ to direct summand $\HEu_i$.
Let $\pi':A\to B(\HEu')$ and $\sigma':B\to B(\HEu')$ be the maps defined by
\begin{align*}
\pi'(a)_{ij}&=\pi(v_i^*av_j) \\
\sigma'(b)_{ij}&=\sigma(v_i^*bv_j).
\end{align*}
Then $\pi'$ and $\sigma'$ are unital $*$--homomorphisms whose restrictions to $D$ agree.
Hence there is a $*$--homomorphism $\psi':\Afr\to B(\HEu')$ with $\psi'\circ\iota_A=\pi'$ and
$\psi'\circ\iota_B=\sigma'$.
Identitying $\HEu$ with $\pi'(p)\HEu'$ and
letting $\psi=\psi'\restrict_{p\Afr p}$ yields the desired $*$--homomorphism.
\end{proof}

\begin{thm}
\label{thm:rfd}
Let $n,n'\in\Nats$.
Suppose $D$ is a finite dimensional C$^*$--algebra and $\gamma:D\to M_n$ and
$\gamma':D\to M_{n'}$ are unital, injective $*$--homomorphisms such that
$\tr_n\circ\gamma=\tr_{n'}\circ\gamma'$, where $\tr_k$ denotes the tracial state on
the $k\times k$ complex matrices $M_k$.
Let
\begin{equation}
\label{eq:Afp}
\Afr= M_n*_DM_{n'}
\end{equation}
be the full free product C$^*$--algebra with amalgamation over $D$,
with respect to the inclusions $\gamma$ and $\gamma'$.
Then $\Afr$ is a residually finite dimensional C$^*$--algebra.
\end{thm}
\begin{proof}
Compressing $\Afr$ with an appropriate projection $p\in D$ and
using Lemmas~\ref{lem:fpcorner} and~\ref{lem:pAp}, we may without loss of generality assume $D$ is abelian.
Let $p_1,p_2,\ldots,p_m$ denote the minimal projections of $D$ and let $t_i=\tr_n\circ\gamma(p_i)$.
Note that $nt_i$ and $n't_i$ are both integers.
Let $\iota:M_n\to\Afr$ and $\iota':M_{n'}\to\Afr$ be the inclusions arising
from the free product construction~\eqref{eq:Afp}.
Let $\pi:\Afr\to B(\HEu)$ be a faithful, unital $*$--homomorphism.
Let $\alpha=\pi\circ\iota$ and $\beta=\pi\circ\iota'$ be the resulting representations of $M_n$
and $M_{n'}$ on the Hilbert space $\HEu$.
Let $F_1\subseteq F_2\subseteq\cdots$ and  $G_1\subseteq G_2\subseteq\cdots$
be chains of finite dimensional subspaces of $\HEu$, each of whose unions is dense in $\HEu$,
and such that each $F_k$ is $\alpha(M_n)$--invariant and each $G_k$ is $\beta(M_{n'})$--invariant.
Let $E_k=F_k+G_k$.

For each $k$, we will now find a finite dimensional subspace $E'_k$ of $\HEu$ that is orthogonal to $E_k$ and
such that,
letting $H_k=E_k+E'_k$, there are $*$--representations $\alphat:M_n\to B(H_k)$ and
$\betat:M_{n'}\to B(H_k)$ such that
\begin{align}
\alphat\circ\gamma&=\betat\circ\gamma' \label{eq:agbg} \\
\alphat(\cdot)\restrict_{F_k}&=\alpha(\cdot)\restrict_{F_k} \label{eq:alres} \\
\betat(\cdot)\restrict_{G_k}&=\beta(\cdot)\restrict_{G_k}\;. \label{eq:betres}
\end{align}
Note that $\alpha\circ\gamma=\beta\circ\gamma'$ and $E_k$ is $\alpha\circ\gamma(D)$--invariant.
Let $d=\dim(E_k)$, $P_i=\alpha\circ\gamma(p_i)\restrict_{E_k}$ and $r_i=\rank(P_i)$.
Let $E'_k$ be any finite dimensional subspace of $\HEu$ that is orthogonal to $E_k$
and having dimension $d'=\dim(E'_k)>0$ so that $d+d'$ is divisible by both $n$ and $n'$
and such that $t_i(d+d')\ge r_i$ for all $i$.
Since each $t_i(d+d')$ is an integer, we can find projections $Q_1,\ldots,Q_m\in B(E_k')$
such that $Q_1+\cdots+Q_m=1$ and, letting $r_i'=\rank(Q_i)$, $r_i+r_i'=t_i(d+d')$.
We will define a representation $\alphat:M_n\to B(H_k)$ satisfying $\alphat\circ\gamma(p_i)=P_i+Q_i$
by setting
\[
\alphat(x)=\alpha(x)\restrict_{F_k}+\alpha'(x)
\]
where $\alpha'$ is any representation of $M_n$ on $H_k\ominus F_k$ such that
$\alpha'\circ\gamma(p_i)=(P_i+Q_i)\restrict_{H_k\ominus F_k}$.
Such a representation $\alpha'$ can always be found because, letting $\dim(F_k)=na$, we have
\[
\dim(H_k\ominus F_k)=d+d'-na=na'
\]
for some positive integer $a'$ and
\[
\rank((P_i+Q_i)\restrict_{H_k\ominus F_k})=r_i+r_i'-t_ina=t_ina'.
\]
By construction, equation~\eqref{eq:alres} holds.
In a similar way, we can find a representation $\betat:M_{n'}\to B(H_k)$
such that $\betat\circ\gamma'(p_i)=P_i+Q_i$ and equation~\eqref{eq:betres} holds.

For a given $k$, having found $\alphat$ and $\betat$, let $\pi_k:\Afr\to B(H_k)$ be the $*$--representation
such that $\pi_k\circ\iota=\alphat$ and $\pi_k\circ\iota'=\betat$.
Let $\Afr_0$ be the dense subalgbra of $\Afr$ that is algebraically generated by $\iota(M_n)\cup\iota'(M_{n'})$.
We will show that for any given $x\in\Afr_0$ and any $\eps>0$ there is $k\in\Nats$ such that
\begin{equation}
\label{eq:pix}
\|\pi_k(x)\|\ge\|x\|-\eps\;.
\end{equation}
This will suffice to show that $\Afr$ is residually finite dimensional.
Let $\xi\in\HEu$ be a unit vector such that $\|\pi(x)\xi\|\ge\|x\|-\eps/2$.
We may write $x=w_1+w_2+\cdots+w_M$ as the sum of finitely many words $w_i$ in $\iota(M_n)$ and $\iota'(M_{n'})$.
We will show that for every $i$ there is $k(i)$ such that if $k\ge k(i)$ then
\begin{equation}
\label{eq:wi}
\|\pi_k(w_i)\xi-\pi(w_i)\xi\|<\eps/2M\;.
\end{equation}
Taking $k\ge\max_{1\le i\le M}k(i)$, this will imply
$\|\pi_k(x)\xi-\pi(x)\xi\|<\eps/2$, which will yield~\eqref{eq:pix}.
To show~\eqref{eq:wi} for fixed $i$, write
\[
w_i=a_\ell a_{\ell-1}\cdots a_2a_1
\]
for some $\ell\in\Nats$ and $a_1,\ldots,a_\ell\in\iota(M_n)\cup\iota'(M_{n'})$.
Let $\xi_0=\xi$ and $\xi_j=\pi(a_j)\xi_{j-1}$, ($1\le j\le\ell$).
Let $N=\max_{1\le j\le\ell}\|a_j\|$.
Choosing $k$ large enough ensures that
\[
\max\big(\dist(\xi_{j-1},F_k),\dist(\xi_{j-1},G_k)\big)<\eps/(4\ell MN^{\ell-j+1})
\]
for all $j\in\{1,\ldots,\ell\}$.
Then the estimates
\begin{align*}
\|\pi(w_i)\xi-\pi_k(w_i)\xi\|&\le\sum_{j=1}^\ell\|\pi_k(a_\ell\cdots a_{j+1})\|\,\|\pi(a_j)\xi_{j-1}-\pi_k(a_j)\xi_{j-1}\| \\
\|\pi(a_j)\eta-\pi_k(a_j)\eta\|&\le2\|a_j\|\max\big(\dist(\eta,F_k),\dist(\eta,G_k)\big)
\end{align*}
give~\eqref{eq:wi}.
\end{proof}

The next example shows that the
condition in Theorem~\ref{thm:rfd} that the traces on $M_n$ and $M_{n'}$ restrict to the same trace on $D$
cannot be dropped.
\begin{example}\rm
\label{ex:notrfd}
Let $D=\Cpx\oplus\Cpx$ and let $p_1$ and $p_2$ be the minimal projections of $D$.
Let $\gamma:D\to M_2$ identify $D$ with the diagonal subalgebra of $M_2$ and let $\gamma':D\to M_3$
be the unital $*$--homomorphism sending $p_1$ to a minimal projection in $M_3$ and thus sending $p_2$
the the sum of two minimal projections in $M_3$.
Let
\[
\Afr=M_2*_DM_3
\]
be the correspoinding full amalgamated free product of C$^*$--algebras.
Since we can easily find representations of $M_2$ and $M_3$ on an infinite dimensional Hilbert space that agree
on the embedded copies of $D$, $\Afr$ is a nonzero C$^*$--algebra.
However, $p_2$ is the sum of two projections in $\Afr$, each equivalent to $p_2$ itself.
Thus any finite dimensional representation of $\Afr$ sends $p_2$ to zero, and therefore must send everything to zero.
We conclude that $\Afr$ has no nonzero finite dimensional representations, and is certainly not residually finite dimensional.
\end{example}

\section{Amalgamated Free Products and Interpolated Free Group Factors} 

In this section we will show that a free product of 
two copies of the hyperfinite II$_{1}$--factor with amalgamation
over a (possibly infinite) direct sum of matrix algebras can be identified with an
interpolated free group factor.
We will let $\delta_{0}$ denote 
Voiculescu's (modified) free entropy dimension (see \cite{dvv:partIII}).
By recent 
work of Kenley Jung~\cite{jung:hyperfinite}
it is now known that $\delta_{0}$ is an
invariant of hyperfinite von Neumann algebras.  Hence, if $B$ is a 
hyperfinite von Neumann algebra with normal tracial state $\tau$, we 
will let $\delta_{0}(B,\tau)$ denote the modified free entropy
dimension of $B$ with respect to $\tau$. 

\begin{thm}\label{thm:RBR}
Let $R$ be the hyperfinite II$_1$--factor with tracial state $\tau$
and let $B\subseteq R$ be a unital W$^*$--subalgebra.
Assume $B$ is type I and has atomic center.
Let
\begin{equation}\label{eq:MRBR}
\Mcal=R*_BR
\end{equation}
be the amalgamated free product of von Neumann algebras, taken with respect to the
$\tau$--preserving conditional expectation $E_B:R\to B$.
Then $M\cong L(\Fb_s)$, where $s=2-\delta_0(B,\tau\restrict_B)$.
\end{thm}    
\begin{proof}
Let $R^{(1)}$ and $R^{(2)}$ denote the copies of $R$ in the amalgamated free product~\eqref{eq:MRBR},
and let $E^{(1)}$ and $E^{(2)}$ denote the respective copies of the conditional expectation $E_B$.
We more formally write
\[
(\Mcal,E)=(R^{(1)},E^{(1)})*_B(R^{(2)},E^{(2)})\;,
\]
thus denoting by $E:\Mcal\to B$ the free product conditional expectation.
We will denote also by $\tau$ the normal tracial state $\tau\restrict_B\circ E$ on $\Mcal$.

The algebra $B$ is a direct sum of (possibly infinitely many) matrix algebras over the complex numbers,
$B=\oplus_{i\in I}M_{n(i)}$.
We will first argue that it will suffice to prove the theorem in the case that $B$ is abelian.
Let $p_i$ be a minimal projection of the $i$th summand, $M_{n(i)}$, and let $p=\sum_{i\in I}p_i$.
Let $\lambda_i=\tau(p_i)$ and $\lambda=\tau(p)$.
Then
\begin{alignat}{2}
s=&\,\delta_0(B,\tau\restrict_B)&&=1-\sum_{i\in I}\lambda_i^2 \label{eq:s} \\
s':=&\,\delta_0(pBp,\lambda^{-1}\tau\restrict_{pBp})&&=1-\sum_{i\in I}(\lambda_i/\lambda)^2\;. \notag
\end{alignat}
Clearly, $p\Mcal p$ is generated by $pR^{(1)}p$ together with $pR^{(2)}p$, and these are free with respect to
the conditional expectation $E\restrict_{p\Mcal p}$.
Thus $p\Mcal p\cong(pRp)*_{pBp}(pRp)$.
Since $pBp$ is abelian, the theorem in this case will yield $p\Mcal p\cong L(\Fb_{s'})$,
and then the rescaling formula for interpolated free group factors (see~\cite{D:interp} and~\cite{Ra94}),
gives $\Mcal\cong L(\Fb_s)$.

Therefore, we may without loss of generality assume $B$ is abelian, and we will denote by $(p_i)_{i\in I}$
the minimal projections  of $B$, with $\lambda_i=\tau(p_i)$.
We will take $I=\{1,2,\ldots,n\}$ if $I$ is finite, and $I=\Nats$ if $I$ is infinite.
We will also assume $\lambda_1=\max_{i\in I}\lambda_i$.
Given $k\in I$, let
\begin{align*}
p[1,k]&=p_1+\cdots+p_k \\
\lambda[1,k]&=\lambda_1+\cdots+\lambda_k \\
\Mcal[1,k]&=W^*(p[1,k]R^{(1)}p[1,k]\cup p[1,k]R^{(2)}p[1,k])\;.
\end{align*}
Then $\Mcal[1,1]\cong R*R\cong L(\Fb_2)$,
where the second isomorphism is from~\cite{D:interp}.
Fix $k$ such that $k,k+1\in I$.
Let $v\in R^{(1)}$ and $w\in R^{(2)}$ be such that $v^*v=p_{k+1}=w^*w$, $vv^*\le p_1$ and $ww^*\le p_1$.
Then $\Mcal[1,k+1]=W^*(\Mcal[1,k]\cup\{v,w\})$.
Since $\Mcal[1,1]$ is a  factor, there is a partial isometry $u\in\Mcal[1,1]$ such that $u^*u=vv^*$ and $uu^*=ww^*$.
Let $q=vv^*$ and $a=vw^*u$.
Then
\begin{equation}\label{eq:kp1}
\Mcal[1,k+1]=W^*(\Mcal[1,k]\cup\{v,a\})\;,
\end{equation}
so
\begin{equation}\label{eq:qkp1}
q\Mcal[1,k+1]q=W^*(q\Mcal[1,k]q\cup\{a\})\;.
\end{equation}
We will show that, with respect to the trace $\lambda_{k+1}^{-1}\tau\restrict_{q\Mcal q}$, $a$ is a Haar unitary element
in $q\Mcal q$ and the pair
\begin{equation}\label{eq:qMqa}
q\Mcal[1,k]q,\quad\{a\}
\end{equation}
is $*$--free.

For a subalgebra $A\subseteq\Mcal$, let $A\oup=A\cap\ker E$.
We will also use the notation
\[
\Lambdao(S_1,\ldots,S_p)=
\{b_1b_2\cdots b_n\mid n\in\Nats,\,b_j\in S_{i_j},\,i_1\ne i_2,\,i_2\ne i_3,\ldots,i_{n-1}\ne i_n\}
\]
for subsets $S_1,\ldots,S_p$ of an algebra.
Clearly, $a$ is a unitary element of $q\Mcal q$.
In order to prove that $a$ is a Haar unitary and that the pair~\eqref{eq:qMqa} is $*$--free,
it will suffice to show $\Theta\subseteq\ker\tau$, where $\Theta$ is the set of all words in
\[
\Lambdao((q\Mcal[1,k]q)\oup,\{a^k\mid k\in\Nats\}\cup\{(a^*)^k\mid k\in\Nats\})
\]
having at least one occurrence of an $a^k$ or an $(a^*)^k$.
Let $q'=ww^*$.
By rewriting $a=vw^*u$ and $a^*=u^*wv^*$, we find $\Theta\subseteq\Phi$, where $\Phi$
is the set of all words $x=x_1x_2\cdots x_n$ belonging to
\[
\Lambdao((q\Mcal[1,k]q)\oup\cup q\Mcal[1,k]q'\cup q'\Mcal[1,k]q\cup(q'\Mcal[1,k]q')\oup,\{vw^*,wv^*\})
\]
such that
\renewcommand{\labelenumi}{$\bullet$}
\begin{enumerate}
\item if $x_{j-1}=wv^*$ and $x_{j+1}=vw^*$, then $x_j\in(q\Mcal[1,k]q)\oup$
\item if $x_{j-1}=vw^*$ and $x_{j+1}=wv^*$, then $x_j\in(q'\Mcal[1,k]q')\oup$
\item there is at least one occurrence of $vw^*$ or $wv^*$.
\end{enumerate}
The subalgebra
\[
B+\lspan\Lambdao((p[1,k]R^{(1)}p[1,k])\oup,(p[1,k]R^{(2)}p[1,k])\oup)
\]
is s.o.--dense in $\Mcal[1,k]$.
Since $q\in R^{(1)}$ and $q\le p_1$ with $p_1$ a minimal projection in $B$,
Kaplansky's density theorem can be used to show that every element of $(q\Mcal[1,k]q)\oup$
is the s.o.--limit of a bounded sequence in
\[
(qR^{(1)}q)\oup+\lspan\big(\Lambdao((p[1,k]R^{(1)}p[1,k])\oup,(p[1,k]R^{(2)}p[1,k])\oup)\backslash(p[1,k]R^{(1)}p[1,k])\oup\big)\;.
\]
Similarly, every element of $(q'\Mcal[1,k]q')\oup$
is the s.o.--limit of a bounded sequence in
\[
(q'R^{(2)}q')\oup+\lspan\big(\Lambdao((p[1,k]R^{(1)}p[1,k])\oup,(p[1,k]R^{(2)}p[1,k])\oup)\backslash(p[1,k]R^{(2)}p[1,k])\oup\big)\;,
\]
and every element of $q\Mcal[1,k]q'\cup q'\Mcal[1,k]q$ is the s.o.--limit of a bounded sequence in
\[
\Cpx p_1+\lspan\Lambdao((p[1,k]R^{(1)}p[1,k])\oup,(p[1,k]R^{(2)}p[1,k])\oup)\;.
\]
Therefore, in order to show $\Phi\subseteq\ker\tau$, it will suffice to show $\Psi\subseteq\ker\tau$,
where $\Psi$ is the set of all words $y=y_1y_2\cdots y_n$ belonging to
\[
\Lambdao\big(\Lambdao((p[1,k]R^{(1)}p[1,k])\oup,(p[1,k]R^{(2)}p[1,k])\oup)\cup\{p_1\},\{vw^*,wv^*\}\big)
\]
such that
\renewcommand{\labelenumi}{$\bullet$}
\begin{enumerate}
\item if $y_{j-1}=wv^*$, $y_{j+1}=vw^*$ and if $y_j\in R^{(1)}$, then $y_j\in(qR^{(1)}q)\oup$
\item if $y_{j-1}=vw^*$, $y_{j+1}=wv^*$ and if $y_j\in R^{(2)}$, then $y_j\in(q'R^{(2)}q')\oup$
\item there is at least one occurrence of $vw^*$ or $wv^*$.
\end{enumerate}
By using the inclusions
\begin{alignat*}{2}
v^*(qR^{(1)}q)\oup v&\subseteq(R^{(1)})\oup\qquad&w^*(q'R^{(2)}q')\oup w&\subseteq(R^{(2)})\oup \\
v^*(p[1,k]R^{(1)}p[1,k])&\subseteq(R^{(1)})\oup\qquad&w^*(p[1,k]R^{(2)}p[1,k])&\subseteq(R^{(2)})\oup \\
(p[1,k]R^{(1)}p[1,k])v&\subseteq(R^{(1)})\oup\qquad&(p[1,k]R^{(2)}p[1,k])w&\subseteq(R^{(2)})\oup\;,
\end{alignat*}
we find $\Psi\subseteq\Lambdao((R^{(1)})\oup,(R^{(2)})\oup)\subseteq\ker E\subseteq\ker\tau$,
where the second inclusion above is by freeness of $R^{(1)}$ and $R^{(2)}$.
This completes the proof that $a$ is Haar unitary and that the pair~\eqref{eq:qMqa} is $*$--free.

We have, therefore, $q\Mcal[1,k+1]q\cong q\Mcal[1,k]q*L(\Ints)$.
Beginning with $\Mcal[1,1]\cong L(\Fb_2)$ and using~\eqref{eq:kp1} and~\eqref{eq:qkp1}, it
follows readily by induction on $k$ that $\Mcal[1,k]$ is an interpolated free group factor with
$p_1\Mcal[1,k]p_1\cong L(\Fb_{s(k)})$, where $s(k)=1+\sum_{j=1}^k(\lambda_j/\lambda_1)^2$,
and that the inclusion $p_1\Mcal[1,k]p_1\hookrightarrow p_1\Mcal[1,k+1]p_1$
is a standard embedding of interpolated free group factors (see~\cite{D:fdim}).
Refering to~\cite{D:fdim} when $I=\Nats$ for the inductive limit of standard embeddings, we conclude
$p_1\Mcal p_1\cong L(\Fb_{s''})$, where $s''=1+\sum_{i\in I}(\lambda_i/\lambda_1)^2$.
The rescaling formula gives $\Mcal\cong L(\Fb_s)$ with $s$ as in~\eqref{eq:s}.
\end{proof}

It will be convenient to state the special case of the theorem above 
which we will need in the next section.
Given integers $k(n)\ge2$ and $\ell(n)\in\{1,2,\ldots,k(n)-1\}$, ($n\in\Nats$),
we consider unital subalgebras $B_n=\Cpx\oplus M_{\ell(n)}\subseteq M_{k(n)}$,
where the summand $M_{\ell(n)}$ of $B_n$ is a corner of $M_{k(n)}$.
We consider the W$^*$--subalgebra $B$ of the hyperfinite II$_1$--factor given by
\[
B=\overline{\otimes}_{n=1}^{\,\infty} B_n
\subseteq\overline{\otimes}_{n=1}^{\,\infty} M_{k(n)}=R\;.
\]
Let $\alpha_n=\ell(n)/k(n)$.
It is straightforward to show that $B$ is type I with atomic center if and only if
$\sum_{n=1}^\infty\alpha_n<\infty$, and that then
\[
B=\bigoplus_{\substack{F\subseteq\Nats \\|F|<\infty}}\underset{\lambda_F}{M_{m(f)}}\;,
\]
where $\lambda_F$ denotes the trace of a minimal projection in the summand $M_{m(F)}$ and
where, using the convention that the empty product is equal to $1$, we have
\[
m(F)=\prod_{n\in F}\ell(n)\,,\quad\lambda_F=\lambda\bigg/\prod_{n\in F}k(n)(1-\alpha_n)\;,
\]
with $\lambda=\prod_{n=1}^\infty(1-\alpha_n)$.
Hence, letting $\tau$ denote the tracial state on $R$, we have
\begin{align*}
\delta_0(B,\tau\restrict_B)&=1-\sum_{\substack{F\subseteq\Nats\\|F|<\infty}}\lambda_F^2 \\
&=1-\lambda^2\prod_{n=1}^\infty\bigg(1+\left(\frac1{k(n)(1-\alpha_n)}\right)^2\bigg) \\
&=1-\prod_{n=1}^\infty\bigg((1-\alpha_n)^2+\frac1{k(n)^2}\bigg)\;.
\end{align*}

\begin{cor}\label{thm:maincor}
In the situation above, assuming $\sum_{n=1}^\infty\alpha_n<\infty$, we have 
$R*_{B} R \cong L(\Fb_{t})$,
where 
\begin{equation}
\label{eq:tval}
t = 1 + \prod_{n = 1}^{\infty} \bigg( (1 - \alpha_{n})^{2} + \frac1{k(n)^2} \bigg)\;.
\end{equation}
Hence, if $p_{1}$ denotes the unit of $M_{\ell(1)}$ (which is 
a projection in $R$ of trace $\alpha_{1}$), then the compression of 
$R*_{B}R$ by $p_{1}$ is isomorphic to  $L(\Fb_s),$ where 
\[
s = 1 + \bigg(\big(\frac1{\alpha_{1}}-1\big)^{2} + \frac1{\ell(1)^2}\bigg) \prod_{n = 
2}^{\infty} \bigg( (1 - \alpha_{n})^{2} + \frac1{k(n)^2} \bigg)\;.
\]
\end{cor}

\section{Construction of Popa Algebras}  

In this section we will construct Popa algebras that are weakly densely embedded in prescribed 
interpolated free group factors.
These factors will arise as in Corollary~\ref{thm:maincor}.
The construction is rather technical and hence we begin this section 
by describing the abstract properties we are after.  

Given $1 < s < \infty$ we will construct an inductive system 
$$A(1) \to A(2) \to A(3) \to \cdots,$$ where each $A(n)$ will be a 
full amalgamated free product of the following form: $$A(n) \cong  
\otimes_{p=1}^{n} M_{k(p)} ({\Cpx}) *_{\otimes_{p=1}^{n} B_{p}} 
\otimes_{p=1}^{n} M_{k(p)} ({\Cpx}).$$  The subalgebras $B_{p} 
\subset M_{k(p)}$ 
will all be of the form $B_p \cong {\Cpx} \oplus M_{\ell(p)}$,
where the summand $M_{\ell(p)}$ is a corner of $M_{k(p)}$.
The sequences $\{ k(p) \}$ and $\{ \ell(p) \}$ will be chosen with care, as we will have to 
arrange the identity
$$s =  1 + \bigg(\big(\frac{1}{\alpha_{1}}-1\big)^{2} + \frac{1}{\ell(1)^{2}} \bigg)
\prod_{n = 2}^{\infty} \bigg( (1 - \alpha_{n})^{2} + \frac1{k(n)^2} \bigg), $$  
where $\alpha_n=\ell(n)/k(n)$.

Letting $\mathfrak{A}$ denote the inductive limit of the sequence $A(n) 
\to A(n+1)$ we will then show, using a von Neumann algebra version of Elliott's intertwining argument,
that $\mathfrak{A}$ has a tracial state 
$\tau$ such that the von Neumann algebra generated by the GNS representation is isomorphic to a (reduced) amalgamated
free product of von Neumann algebras, namely
$$\pi_{\tau} (\mathfrak{A})^{\prime\prime} \cong 
\bar{\otimes}_{p=1}^{\infty} M_{k(p)} *_{\bar{\otimes}_{p=1}^{\infty} B_{p}}  
\bar{\otimes}_{p=1}^{\infty} M_{k(p)}.$$  Finally, we will show that 
$\mathfrak{A}$ is a Popa algebra and that if $p_{1} \in \mathfrak{A}$ denotes 
the unit of $M_{\ell(1)} \subset B_{1} \subset \mathfrak{A}$ then the 
corner $p_{1}\mathfrak{A}p_{1}$ is again a Popa algebra, denoted by $A_{s}$, 
which, by Corollary 
\ref{thm:maincor}, has a GNS representation isomorphic to $L({\Fb}_{s})$.
Since $A_s$ is simple, this GNS representation gives an embedding of $A_s$ as a weakly dense subalgebra of $L(\Fb_s)$.

More formally, what follows is the main theorem of this section.

\begin{thm}
\label{thm:paconstruction}
For every $1 < s < \infty$ there exists an inductive system of C$^*$--algebras
$A(1) \to A(2) \to \cdots$ with inductive limit C$^*$--algebra $\mathfrak{A}$ such 
that:  
\renewcommand{\labelenumi}{(\arabic{enumi})}
\begin{enumerate}
 \item $A(n)$ is the full amalgamated free product C$^*$--algebra
 $$A(n) = \big(\otimes_{p = 1}^{n} M_{k(p)} \big) *_{\otimes_{p = 
 1}^{n} B_{p}}  \big( \otimes_{p = 1}^{n} M_{k(p)} \big), $$
 where $B_{p}=\Cpx\oplus M_{\ell(p)} \subset 
 M_{k(p)}$ is a unital subalgebra with the summand $M_{\ell(p)}$ being
 a corner of $M_{k(p)}$,  for appropriately chosen integers $k(p)$ and $\ell(p)$.
 
 \item Letting $\rho_{m,n}: A(n) \to A(m)$ ($n \leq m$) denote the 
 connecting maps in the inductive system and $p_{m} \in B_{m} =  
 {\Cpx} \oplus  M_{\ell(m)}$ be the unit of the summand $M_{\ell(m)}$, we have that 
 each $\rho_{m,n}$ is injective, 
 $[p_{m}, \rho_{m,n}(x)] = 0$ for all $x \in A(n)$,  
 $p_{m} \rho_{m,n}(x) \in M_{\ell(m)} \subset B_{m}$ for all $x \in A(n)$ 
 and, finally, for each $x \in A(n)$, $\| x \| = \lim_{m\to \infty} \| 
 p_{m} \rho_{m,n}(x) \|$.

 \item Letting $\alpha_n=\ell(n)/k(n)$, we have
 $\alpha_{n} <  2^{-(n+p)}$, for some fixed integer $p$ and 
 $$s =  1 + \bigg(\big(\frac{1}{\alpha_{1}}-1\big)^{2} + \frac{1}{\ell(1)^{2}} \bigg)
 \prod_{n = 2}^{\infty} \bigg( (1 - \alpha_{n})^{2} + \frac1{k(n)^2}
 \bigg) .$$

 \item $\mathfrak{A}$ is a Popa algebra generated by four self adjoint 
 elements and there exists a tracial 
 state $\tau$ on $\mathfrak{A}$ whose GNS representation generates a von Neumann algebra
\begin{equation}\label{eq:amalgfpvN}
\pi_{\tau} (\mathfrak{A})^{\prime\prime} \cong 
\bar{\otimes}_{p=1}^{\infty} M_{k(p)} *_{\bar{\otimes}_{p=1}^{\infty} B_{p}}  
\bar{\otimes}_{p=1}^{\infty} M_{k(p)},
\end{equation}
isomorphic to the indicated amalgamated free product of von Neumann algebras, taken with respect to the
trace--preserving conditional expectations.

 \item Letting $A_{s} = p_{1}\mathfrak{A}p_{1}$ we have that $A_{s}$ 
 is a Popa algebra such that $\pi_{\tau\restrict_{A_{s}}} 
 (A_{s})^{\prime\prime} \cong L(\Fb_{s})$. 
\end{enumerate}
    
\end{thm}

\begin{proof}
Fix $1 < s < \infty$.
Choose a rational number $\alpha_{1} = p(1)/q(1) < 1$, with $p(1), q(1) \in 
{\Nats}$, such that $$s + 1/4 < 1 + (1 - \frac{1}{\alpha_{1}})^{2} < 
s + 1/2.$$  Choose a natural number 
$j(1)$ so large that $$s + 1/4 < 1 + \bigg((1 - \frac{1}{\alpha_{1}})^{2} 
+ \frac{1}{(p(1)j(1))^{2}} \bigg)    < s + 1/2.$$ Then let $k(1) = 
q(1)j(1)$, $\ell(1) = p(1)j(1)$ and $B_{1} =  {\Cpx} \oplus 
M_{\ell(1)}  \cong  {\Cpx} \oplus M_{p(1)} \otimes 
M_{j(1)} \subset M_{q(1)} \otimes M_{j(1)} = M_{k(1)}.$ (We choose 
this embedding to be unital and such that the unit of $M_{\ell(1)}$ is a 
projection in $M_{k(1)}$ of trace $\alpha_{1} = p(1)/q(1)$.)  

Let $A(1)$ denote the full amalgamated free product $M_{k(1)} *_{B_{1}} 
M_{k(1)}$.  By Theorem \ref{thm:rfd},  $A(1)$ is a residually finite dimensional 
C$^{*}$-algebra.  Fix a sequence $\{ a^{(1)}_{m} \}$ which is dense 
in the unit ball of $A(1)$.  Now let $\pit_{1} : A(1) \to M_{t(1)} 
({\Cpx})$ be a unital *-homomorphism such that $\| \pit_{1} 
(a^{(1)}_{1})\| \geq (1/2)\| a^{(1)}_{1} \|$.  Choose a rational number 
$\alpha_{2} = p(2)/q(2) < 1$ such that
$$s + 1/8 < 1 + 
\bigg((1 - \frac{1}{\alpha_{1}})^{2} 
+ \frac{1}{\ell(1)^{2}} \bigg) (1 - \alpha_{2})^{2}  < s + 1/4.$$  
Choose a natural number $j(2)$ so large that $$ s + 1/8 < 1 + 
\bigg((1 - \frac{1}{\alpha_{1}})^{2} 
+ \frac{1}{\ell(1)^{2}} \bigg) \bigg( (1 - \alpha_{2})^{2} + 
(\frac{\alpha_{2}}{j(2)p(2)t(1)})^{2}\bigg)  < s + 1/4.$$
Define $\ell(2) = j(2)p(2)t(1)$, $k(2) = j(2) q(2)t(1)$ and 
$B_{2} =  {\Cpx} \oplus 
M_{\ell(2)}  \cong  {\Cpx} \oplus M_{p(2)} \otimes M_{t(1)} \otimes 
M_{j(2)} \subset M_{q(2)} \otimes M_{t(1)}\otimes M_{j(2)} = M_{k(2)}.$ 
(We choose 
this embedding to be unital and such that the unit of $M_{\ell(2)}$ is a 
projection in $M_{k(2)}$ of trace $\alpha_{2}$.)

Let $A(2)$ be the full amalgamated free product  $M_{k(1)} \otimes 
M_{k(2)} *_{B_{1}\otimes B_{2}}  M_{k(1)}\otimes M_{k(2)}$.  We define 
a $*$-homomorphism $\rho_{2,1}: A(1) \to A(2)$ by the mapping $$x 
\mapsto 
(1 - p_{2}) \sigma_{1}(x) (1 - p_{2}) \oplus \pi_{1} (x),$$ where 
$\sigma_{1}: A(1) \to A(2)$ is the {\em canonical} unital 
*-homomorphism (i.e.\ the natural identification $ M_{k(1)} *_{B_{1}}  
M_{k(1)} \cong M_{k(1)} \otimes 
1_{k(2)} *_{B_{1}\otimes 1_{k(2)}}  M_{k(1)}\otimes 1_{k(2)}
\subset M_{k(1)} \otimes 
M_{k(2)} *_{B_{1}\otimes B_{2}}  M_{k(1)}\otimes M_{k(2)}$), 
$p_{2}$ is the unit of $M_{\ell(2)} \subset A(2)$ and 
$\pi_{1} (x)=\pit_1(x)\otimes1\otimes1 \in M_{t(1)}\otimes M_{p(2)}\otimes 
M_{j(2)} = 
M_{\ell(2)}$.
Note that $\rho_{2,1}$ is a unital *-monomorphism (it is 
not hard to see that $\sigma_{1}$ is injective) whose image commutes 
with $p_{2}$ and such that $p_{2}\rho_{2,1}(x) \in M_{\ell(2)}$ for all 
$x \in A(1)$.  

Note also that if $\tau_{2}$ denotes the canonical free product 
tracial state on $A(2)$ then $\tau_{2} (\rho_{2,1}(1 - p_{1}) (1 - 
p_{2})) = (1 - \alpha_{1})(1 - \alpha_{2})$.  (This observation will 
be relevant later on.)  

The remainder of the 
construction recursively follows a similar pattern. Indeed, we claim 
that following the pattern above one can choose rational numbers $\alpha_n\in(0,1)$ and
integers $k(n)$ and $\ell(n)$ and one can construct algebras $A(n)$ and 
injective $*$-homomorphisms $\rho_{n+1,n} : A(n) \hookrightarrow 
A(n+1)$ with all of the following properties.

Firstly, we have
$$A(n) = \big(\otimes_{m = 1}^{n} M_{k(m)} \big) *_{\otimes_{m = 
    1}^{n} B_{m}}  \big( \otimes_{m = 1}^{n} M_{k(m)} \big),$$ where 
       \begin{enumerate}
	   \item[(a)] $B_{m} \cong {\Cpx} \oplus M_{\ell(m)}$;
	   
	   \item[(b)] the inclusion $B_{m} \subset M_{k(m)}$ is unital and 
	   is such that the unit of $M_{\ell(m)}$ is a projection in $M_{k(m)}$ 
	   of trace $\alpha_{m}$;
	   
	   \item[(c)] each $\alpha_{m}$  is chosen so that the inequality
	   $$ s + 2^{-(m+1)} < 1 + \bigg((1 - \frac{1}{\alpha_{1}})^{2} 
           + \frac{1}{\ell(1)^{2}} \bigg)\bigg(\prod_{i=2}^{m-1}   
	   \big((1 - \alpha_{i})^{2} + \frac1{k(i)^2}\big)\bigg) 
	   (1 - \alpha_{m})^{2} < s + 2^{-m} $$ 
           holds
	   
           \item[(d)] the positive integers $\ell(m)$ and $k(m)$ (which are chosen after 
           $\alpha_{m}$) have the property that $\ell(m)/k(m) = \alpha_{m}$ 
           and, moreover, satisfy the inequalities $$ s + 
           2^{-(m+1)} < 1 + \bigg((1 - \frac{1}{\alpha_{1}})^{2} 
           + \frac{1}{\ell(1)^{2}} \bigg)\prod_{i=2}^{m} \bigg((1 - \alpha_{i})^{2} + 
           \frac1{k(i)^2}\bigg)  < s + 2^{-m}.$$ 
       \end{enumerate}
Moreover, the connecting maps $\rho_{n+1,n} : A(n) \to A(n+1)$ are all 
    of the form $$x \mapsto (1 - p_{n+1}) \sigma_{n}(x) (1 - p_{n+1}) \oplus 
    \pi_{n} (x),$$ where 
        \begin{enumerate}
            \item[(e)] $\sigma_{n} : A(n) \to A(n+1)$ is the {\em 
            canonical} unital $*$-monomorphism,
            i.e.\ the natural identification
            \begin{align*}
             \otimes_{p=1}^n&M_{k(p)}*_{\otimes_{p=1}^nB_{p}}  
            \otimes_{p=1}^nM_{k(p)} \\
            &\cong
            \big((\otimes_{p=1}^nM_{k(p)})\otimes1_{k(n+1)}\big)*_{(\otimes_{p=1}^nB_{p})\otimes1_{k(n+1)}}  
            \big((\otimes_{p=1}^nM_{k(p)})\otimes1_{k(n+1)}\big) \\
            &\subset
            \otimes_{p=1}^{n+1}M_{k(p)} *_{\otimes_{p=1}^{n+1}B_{p}}  
            \otimes_{p=1}^{n+1}M_{k(p)}
            \end{align*}
            of full amalgamated free products;
	    
	    \item[(f)] $p_{n+1}$ is the unit of $M_{\ell(n+1)} \subset B_{n+1} 
	    \subset A(n+1)$; 
	    
	    \item[(g)] $\pi_{n} : A(n) \to M_{\ell(n+1)}$ 
	    is a unital $*$-monomorphism with the property   
	    that $\| \pi_{n} (\rho_{n,t}(a_{s}^{(t)})) \| 
	    \geq (1 - 2^{-n}) \| 
	    a_{s}^{(t)} \|,$ for all $1 \leq s,t \leq n-1$, where 
	    $\{a_{s}^{(t)} \}_{s \in {\Nats}}$ is a sequence which is 
	    dense in the unit ball of $A(t)$ and $\rho_{n,t} = 
	    \rho_{n,n-1}\circ \cdots \circ \rho_{t+1,t}: A(t) \to A(n)$.
        \end{enumerate}  

To perform the induction step, assume that algebras $A(1), \ldots, A(n-1)$ and connecting maps 
$\rho_{i+1,i} : A(i) \to A(i+1)$, $i \in \{1,\ldots,n-2 \}$, have been 
constructed with all the properties above.  We will construct $A(n)$ 
and the appropriate connecting map $\rho_{n,n-1}$.  

First, since by Theorem 
\ref{thm:rfd} $A(n-1)$ is residually finite dimensional, 
we can find a finite dimensional representation  $\tilde{\pi}_{n-1} : 
A(n-1) \to M_{t(n-1)} 
({\Cpx})$ such that $\| \tilde{\pi}_{n-1} 
(\rho_{n-1,t} (a^{(t)}_{s}))\| \geq (1 - 2^{-n})\| a^{(t)}_{s} \|$ 
for all $1 \leq s,t \leq n-1$.  Choose a rational number 
$\alpha_{n} = p(n)/q(n) < 1$ satisfying the inequalities in (c) above.
Choose a natural number $j(n)$ so large that when one defines 
$\ell(n) = j(n)p(n)t(n-1)$ and $k(n) = j(n)q(n)t(n-1)$ one gets the 
inequalities in part (d).  The rest of the construction proceeds along the lines of the case
$n=2$ treated above, in order that the desired properties hold.

Now let $\mathfrak{A}$ denote the inductive limit of the inductive 
system $\{A(n), \rho_{m,n}\}$.  It is clear from the construction that 
we have satisfied all of the statements in parts (1), (2) and (3) 
of the theorem except for the inequalities claimed in part (3).  
Hence parts (1), (2) and (3) will be complete as soon as we prove 
the following inequality: $$\alpha_{m} 
<  2^{-(m+p)},$$ where $p$ is some integer.
Let $p\in\Ints$ be such that $s -1 > 
2^{(p+1)}$ (since $s > 1$, such a $p$ exists).  
Letting $$\gamma_{m-1} =  \bigg((1 - 
\frac{1}{\alpha_{1}})^{2} + \frac{1}{\ell(1)^{2}} \bigg)\prod_{i=2}^{m-1}  
\bigg((1 - \alpha_{i})^{2} + \frac1{k(i)^2}\bigg)$$ we have 
the inequalities $s + 2^{-m} < 1 + \gamma_{m-1} < s + 2^{-(m-1)}$ 
and $s + 2^{-(m+1)} < 1 + \gamma_{m-1} (1 - \alpha_{m})^{2} < s + 
2^{-m}$.  It follows that $\gamma_{m-1}[1 - (1 - \alpha_{m})^{2}] 
< 2^{-(m-1)} - 2^{-(m+1)} < 2^{-(m-1)}$. Since $\gamma_{m-1} > s - 1 > 
2^{p}$, for some integer $p$ it follows that $$\alpha_{m}(2 - 
\alpha_{m}) < 2^{-(m+p-1)}.$$  Since $\alpha_{m} < 1$ by construction, 
we replace $p$ by $p-1$ and the inequality claimed above is now immediate.
We have thus shown that parts~(1), (2) and~(3) of the theorem hold.

To prove parts (4) and (5) in the statement of the 
theorem, let us first observe that both $\mathfrak{A} = \varinjlim A(n)$ and 
$A_{s} = p_{1} \mathfrak{A} p_{1}$ are Popa algebras.  It is clear 
(from part (2) of the theorem) 
that $\mathfrak{A}$ is a unital C$^{*}$-algebra which satisfies the 
finite dimensional approximation property which defines Popa 
algebras.  Hence the only question is whether or not $\mathfrak{A}$ 
is simple.  However, this also follows from part (2) of the theorem.  
Indeed, if $I \subset \mathfrak{A}$ is a non-zero ideal then we can 
find a non-zero element $x \in A(n) \cap I$ (identifying $A(n)$ with 
its image in $\mathfrak{A}$ and taking a sufficiently large 
$n$).  Choosing $m$ large enough that $0 \neq p_{m}\rho_{m,n}(x) \in 
M_{\ell(m)} \subset B_{m}$ we see that the ideal (of $A(m)$) generated by 
$\rho_{m,n}(x)$ is all of $A(m)$ since $B_{m}$ is contained in a 
unital matrix subalgebra of $A(m)$.  It follows that $I = 
\mathfrak{A}$ and hence $\mathfrak{A}$ is simple.  

To see that $A_{s}$ is also a Popa algebra, we only need to prove 
that $A_{s}$ has the right finite dimensional approximation property 
since $A_{s}$ is clearly unital and simple (being a corner of the 
simple C$^{*}$-algebra $\mathfrak{A}$).  Note that $A_{s} = \varinjlim 
\rho_{n,1}(p_{1}) A(n) \rho_{n,1}(p_{1})$ so to prove Popa's 
approximation property for 
$A_{s}$ it suffices to consider a finite set $\mathfrak{F} 
\subset \rho_{n,1}(p_{1}) A(n) \rho_{n,1}(p_{1})$ for some $n$.  
Note that since $p_{m}\rho_{m,n}(x) \in M_{\ell(m)} \subset B_{m}$ for 
all $x \in A(n)$ ($n < m$) it 
follows  that $\rho_{m,1}(p_{1}) p_{m} B_{m} p_{m}
\rho_{m,1}(p_{1}) \subset M_{\ell(m)}$ is a 
non-zero finite dimensional subalgebra of $A_{s}$ (again, 
identifying $A(m)$ with its image in $\mathfrak{A}$) with the property 
that its unit commutes with $\rho_{m,n}\big(\rho_{n,1}(p_{1}) A(n) 
\rho_{n,1}(p_{1})\big)$ for all $n < m$ and compressing by this unit 
(i.e.\ $\rho_{m,1}(p_{1})p_{m}$) maps 
any element of $\rho_{m,n}\big(\rho_{n,1}(p_{1}) A(n) 
\rho_{n,1}(p_{1})\big)$ into this 
finite dimensional subalgebra.  This is stronger than the 
approximation property defining Popa algebras and hence we see that 
$A_{s}$ is also a Popa algebra.

We also claimed in (4) that $\mathfrak{A}$ is generated by four 
self adjoint elements.  To see this, it suffices to show that 
$\mathfrak{A}$ is generated by two UHF algebras since UHF algebras 
are generated by two self--adjoint elements.   (This last statement
is well known to the experts.
A proof can be given as follows.
Since $C(X)$, the continuous functions on 
the Cantor set $X$, is generated by a single self--adjoint element, it will suffice to show that 
any UHF algebra is generated by two copies of $C(X)$.  If ${\mathcal 
U} = \otimes_{i \in {\Nats}} M_{m(i)}$ then it is not hard to 
see that ${\mathcal U} = C^{*}(C(X)_{1}, C(X)_{2})$ where $C(X)_{1} 
\cong C(X)_{2} \cong C(X)$ and $C(X)_{1}$ is generated by the minimal 
projections in the matrix algebras $M_{m(i)}$, while $C(X)_{2} = 
C^{*} (\{ u_{i} : u_{i} \in M_{m(i)} \})$, where the $u_{i}$'s are 
cyclic permutation matrices.) 
However, since the connecting maps $A(n) \to A(n+1)$ map the matrices 
which generate $A(n)$  into the matrices which generate $A(n+1)$,  
it is easy to see that $\mathfrak{A}$ contains two copies of the UHF 
algebra $\otimes_{n=1}^{\infty} M_{k(n)}$ and that these two UHF 
algebras generate $\mathfrak{A}$.  

Finally we are ready to tackle the problem of the GNS representation.
The remainder of the proof is very similar to the proof of Theorem 
2.5.1  in \cite{brown:AFDtraces}.  To ease 
notation we will, for the rest of the proof, identify each of the 
algebras $A(n)$ with its canonical image in $\mathfrak{A}$ (recall 
that all the connecting maps are injective and hence we have 
canonical embeddings $\Phi_{n}: A(n) \hookrightarrow \mathfrak{A}$).  
Hence if $x \in 
A(n)$, $y \in A(m)$ and $n < m$ then when we write $xy$ we really 
mean $\Phi_{m}(\rho_{m,n}(x) y)$.

First we need to define the correct trace on $\mathfrak{A}$.  Let 
$\tau_{n}$ be the canonical free product tracial state on $A(n)$.
  (By the {\em canonical free product tracial state}
  we mean $\tau_{n} = \tr_{k(1)\cdots k(n)} \circ \Et_{n} \circ 
  \lambda_{n}$, where
  $\tr_{k(1)\cdots k(n)}$ is the unique tracial state on 
  $\otimes_{m=1}^{n} M_{k(m)}$,
  $\lambda_{n} : A(n) \to \big( \otimes_{m=1}^{n} 
  M_{k(m)}, E_{n} \big) *_{\otimes_{m=1}^{n} B_{m}}  \big( \otimes_{m=1}^{n} 
  M_{k(m)}, E_{n} \big)$ is the canonical quotient mapping onto the 
  reduced amalgamated free product C$^*$--algebra with respect to the unique 
  $\tr_{k(1)\cdots k(n)}$-preserving conditional expectation $E_{n} : 
  \otimes_{m=1}^{n} M_{k(m)} \to \otimes_{m=1}^{n} B_{m}$ and
  $\Et_n$ is the free product conditional expectation from the above reduced free product
  C$^*$--algebra to $\otimes_{m=1}^nB_m$.)
Extend each $\tau_{n}$ to a state 
$\gamma_{n} \in S(\mathfrak{A})$ 
on $\mathfrak{A}$.  Let $\tau \in S(\mathfrak{A})$ be any weak-* 
cluster point of the sequence $\{\gamma_{n}\}$.  It is easy to check 
that $\tau$ is actually a tracial state on $\mathfrak{A}$ and this is 
the trace that we will need.

Recall that $p_{n}$ denotes the unit of $M_{\ell(n)} \subset B_{n} 
\subset A(n)$.  Let $q_{n} = p_{n}^{\perp} = 1_{\mathfrak{A}} - 
p_{n}$.  Note that since all the $p_{n}$'s commute the same is true 
for all the $q_{n}$'s.  Hence for each pair of natural numbers $n \le m$ 
we can define a projection $Q^{(n)}_{m} \in A(m)$ by $$Q^{(n)}_{m} = 
q_{n}q_{n+1}\cdots q_{m}.$$  We claim that the projections 
$Q^{(n)}_{m}$ enjoy the following properties: 

\begin{enumerate}
    \item[(h)] $[Q^{(n)}_{m}, x] = 0$ for all $x \in A(n-1)$. 
    
    \item[(i)] For fixed $n$, $\{ Q^{(n)}_{m} \}$ is a decreasing 
    sequence (as $m \to \infty$) of projections in $\mathfrak{A}$. 
    
    \item[(j)] $\tau(Q^{(n)}_{m}) = \tau_{m} (Q^{(n)}_{m}) = 
    (1 - \alpha_{n})(1 - 
    \alpha_{n+1}) \cdots (1 - \alpha_{m})
    \ge\prod_{k=n}^\infty(1-\alpha_k) \to 1$ as $n \to \infty$. 

    \item[(k)]  For all $x \in A(n-1)$ we have the identity 
    $$\frac{\tau_{m}(Q^{(n)}_{m} x)}{\tau(Q^{(n)}_{m})} = \tau_{n-1}(x).$$ 
    
\end{enumerate}     

Statements (h) and (i) above are evident from the construction. 

Statement (j) follows from the observation that $Q^{(n)}_{m}$ sits 
inside a unital matrix subalgebra of $A(m)$ and hence all tracial 
states on 
$A(m)$ agree on $Q^{(n)}_{m}$.
Since $\sum\alpha_k<\infty$, we easily get $\lim_{n\to\infty}\prod_n^\infty(1-\alpha_k)=1$.

Part (k) also follows from the construction since $Q^{(n)}_{m} x = 
Q^{(n)}_{m} \sigma_{m-1} \circ \cdots \circ \sigma_{n-1} (x)$ for 
all $x \in A(n-1)$ and $m > n$.  Hence 
\begin{eqnarray*}
\tau_{m} (Q^{(n)}_{m} x)      
      & = & \tau_{m}(Q_{m}^{(n)} \sigma_{m-1}\cdots\sigma_{n-1}(x)) \\ 
      & = & \tr_{k(1)\cdots k(m)} \big( \Et_{m}( \lambda_{m}(Q^{(n)}_{m} 
            \sigma_{m-1}\cdots\sigma_{n-1}(x))) \big)  \\ 
      &	= & \tr_{k(1)\cdots k(m)} \big( \lambda_{m}(Q^{(n)}_{m}) 
           \Et_{m}(\lambda_{m}(\sigma_{m-1}\cdots\sigma_{n-1}(x))) \big) \\
      & = & \tau_{m}(Q^{(n)}_{m}) 
      \tau_{m}(\sigma_{m-1}\cdots\sigma_{n-1}(x)) \\ 
      & = & \tau( Q^{(n)}_{m})\tau_{n-1}(x),
\end{eqnarray*}
for all $x \in A(n-1)$;
the fourth equality above follows because
\begin{align*}
\lambda_m(Q_m^{(n)})&\in1\otimes\cdots\otimes1\otimes B_n\otimes\cdots\otimes B_m \\
\Et_m(\lambda_m(\sigma_{m-1}\cdots\sigma_{n-1}(x)))&\in B_1\otimes\cdots\otimes B_{n-1}\otimes1\cdots\otimes1.
\end{align*}

Now let $\pi_{\tau} : \mathfrak{A} \to B(H_{\tau})$ be the GNS 
representation and define projections $Q^{(n)} \in 
\pi_{\tau}(\mathfrak{A})^{\prime\prime}$ by $$Q^{(n)} = (s.o.t.) 
\lim_{m \to \infty} \pi_{\tau}(Q^{(n)}_{m}).$$ Note that this limit 
exists by (i) above.  
 
We claim that the projections $Q^{(n)}$ have the following 
properties:  

\begin{enumerate}
    \item[(l)] $Q^{(n)} \in \pi_{\tau}(A(n-1))^{\prime}$.  
    
    \item[(m)] $Q^{(1)} \leq Q^{(2)} \leq \cdots$ and 
    $\tau^{\prime\prime} (Q^{(n)}) \to 1$ as $n \to 
    \infty$, where $\tau^{\prime\prime}$ denotes the vector trace 
    induced by $\tau$.
    
    \item[(n)]  For all $x \in A(n-1)$ we have the identity 
    $$\frac{\tau^{\prime\prime}(Q^{(n)}\pi_\tau(x))}{\tau^{\prime\prime}(Q^{(n)})} 
    = \tau_{n-1}(x).$$ 
    
    \item[(o)] For every $x \in A(n-1)$ we have 
    $Q^{(n)}\pi_{\tau}(x) = Q^{(n+1)} \pi_{\tau} (q_{n}x)=Q^{(n+1)}\pi_\tau(q_n\sigma_{n-1}(x))$.
    Hence 
    we have $Q^{(n)} 
    \pi_{\tau}(A(n-1)) \subset Q^{(n+1)}\pi_{\tau} 
    (A(n))$.
    
\end{enumerate}   

Properties (l), (m) and (o) are immediate from the construction.  To 
see property (n), first note that it follows from property (k) 
that $$\lim_{m \to \infty} \tau_{m} (Q^{(n)}_{m}x) = 
\tau^{\prime\prime} (Q^{(n)}) \tau_{n-1} (x),$$ for all $x \in 
A(n-1)$.  Hence it suffices to show that $| \tau(Q^{(n)}_{m} x) - 
\tau_{m}(Q^{(n)}_{m} x) | \to 0$ as $m\to\infty$ (since 
$\tau^{\prime\prime}(Q^{(n)}\pi_\tau(x)) = \lim_{m} \tau (Q^{(n)}_{m} x)$).  
However, by construction, for any $y \in A(p)$ with $\| 
y \| \leq 1$, we have
$|\tau_{p+1}(\rho_{p+1,p}(y))-\tau_p(y)|\le\alpha_{p+1}$
and therefore
$$| \tau(y) - \tau_{p} (y) | \leq \limsup_{j \to 
\infty} |\tau_{j} (\rho_{j,p}(y)) - \tau_{p}(y) | \leq 
\sum_{j=p+1}^{\infty} \alpha_{j} \to 0,$$ as $p \to \infty$. 

Let
\[
\Dfr=\big(\otimes_{n=1}^\infty M_{k(n)}\big)*_{\otimes_1^\infty B_n}\big(\otimes_{n=1}^\infty M_{k(n)}\big)
\]
be the full amalgmated free product C$^*$--algebra and let $\gamma$ be its canonical free product trace,
i.e.\
$\gamma=(\otimes_1^\infty\tr_{k(n)})\circ\Et\circ\sigma$, where $\otimes_1^\infty\tr_{k(n)}$ is the tracial state
on the UHF algebra $\otimes_{n=1}^\infty M_{k(n)}$, where $\sigma:\Dfr\to D$ is the canonical quotient map onto the reduced
free product C$^*$--algebra
\[
(D,\Et)=\big(\otimes_{n=1}^\infty M_{k(n)},E\big)*_{\otimes_1^\infty B_n}\big(\otimes_{n=1}^\infty M_{k(n)},E\big)
\]
taken with respect to the $\otimes_1^\infty\tr_{k(n)}$--preserving conditional expectation
$E:\otimes_{n=1}^\infty M_{k(n)}\to\otimes_1^\infty B_n$ and where $\Et:D\to\otimes_1^\infty B_n$ is the
free product conditional expectation.
Note that the image of $\Dfr$ under the GNS representation $\pi_\gamma$ associated to $\gamma$ is isomorphic to $D$,
and hence generates a von Neumann algebra
$\pi_\gamma(\Dfr)''$ that is isomorphic to the amalgamated free product of von Neumann algebras
found in the RHS of~\eqref{eq:amalgfpvN}.

We will use Elliott's approximate intertwining argument in the context of von Neumann algebras
(see \cite[Theorem 2.3.1]{brown:AFDtraces}) to prove $\pi_\tau(\Afr)''\cong\pi_\gamma(\Dfr)''$.
The C$^*$--algebra $\Dfr$ is the inductive limit of the system of $*$--monomorphisms
$A(1)\overset{\sigma_1}\rightarrow A(2)\overset{\sigma_2}\rightarrow\cdots$,
where $\sigma_n$ is the canonical embedding described in~(e) above.
Thus we regard $A(n)$ as a unital subalgebra of $\Dfr$.
Let $D_n=\pi_\gamma(A(n))$.
Then $D_n$ is a unital C$^*$--subalgebra of $D_{n+1}$, and the inculsion $D_n\hookrightarrow D_{n+1}$
is given by $\pi_\gamma(x)=\pi_\gamma(\sigma_n(x))$, for $x\in A(n)$.
Moreover, $\cup_{n=1}^\infty D_n$ is weakly dense in $\pi_\gamma(\Dfr)''$.
On the other hand, let $C_n=Q^{(n+1)}\pi_\tau(A(n))$.
By par~(l) above, $C_n$ is a C$^*$--algebra.
By~(o), we have the nonunital inclusion $C_n\subseteq C_{n+1}$.
From~(m), it follows that $\cup_{n=1}^\infty C_n$ is weakly dense in $\pi_\tau(\Afr)''$.
Consider $\phi_n:C_n\to D_n$ and $\psi_n:D_n\to C_{n+1}$ given by
\begin{alignat}{2}
\phi_n&(Q^{(n+1)}\pi_\tau(x))=\pi_\gamma(x)&(x&\in A(n)) \label{eq:phin} \\
\psi_n&(\pi_\gamma(x))=Q^{(n+2)}\pi_\tau(\sigma_n(x))\qquad&(x&\in A(n)). \label{eq:psin}
\end{alignat}
Note that~\eqref{eq:phin} defines a $*$--homomorphism and, in fact, a $*$--isomorphism
from $C_n$ to $D_n$ because, given $x\in A(n)$ and using~(n), we get
$Q^{(n+1)}\pi_\tau(x)=0$ $\Leftrightarrow$
$\tau''(Q^{(n+1)}\pi_\tau(x^*x))=0$ $\Leftrightarrow$
$\tau_n(x^*x)=0$ $\Leftrightarrow$
$\gamma(x^*x)=0$ $\Leftrightarrow$
$\pi_\gamma(x)=0$.
On the other hand, given $x\in A(n)$, we have
\[
\phi_{n+1}(Q^{(n+2)}\pi_\tau(\sigma_n(x)))=\pi_\gamma(\sigma_n(x)),
\]
which is the image of $\pi_\gamma(x)$ under the inclusion $D_n\hookrightarrow D_{n+1}$.
Hence~\eqref{eq:psin} defines a $*$--monomorphism $\psi_n$ and the composition $\phi_{n+1}\circ\psi_n$
is equal to the inclusion $D_n\hookrightarrow D_{n+1}$.
Finally, we have
$\psi_n\circ\phi_n(Q^{(n+1)}\pi_\tau(x))=Q^{(n+2)}\pi_\tau(\sigma_n(x))$,
while by~(o),
$Q^{(n+1)}\pi_\tau(x)=Q^{(n+2)}\pi_\tau(q_{n+1}\sigma_n(x))=\pi_\tau(q_{n+1})Q^{(n+2)}\pi_\tau(\sigma_n(x))$,
so invoking~(j) we get
\begin{align*}
\|\psi_n\circ\phi_n(Q^{(n+1)}\pi_\tau(x))-Q^{(n+1)}\pi_\tau(x)\|_2
&\le\|1-\pi_\tau(q_{n+1}))\|_2\,\|Q^{(n+1)}\pi_\tau(x)\| \\
&=\alpha_{n+1}\|Q^{(n+1)}\pi_\tau(x)\|,
\end{align*}
where the $2$--norm is with respect to $\tau''$,
and $\alpha_{n+1}\to0$ as $n\to\infty$.
Now \cite[Theorem 2.3.1]{brown:AFDtraces} can be used to convert $\phi_n$ and $\psi_n$ into an isomorphism
$\pi_\tau(\Afr)''\cong\pi_\gamma(\Dfr)''$, proving part~(4) of the theorem.

Keeping in mind that $\pi_\gamma(\Dfr)''$ is isomorphic
to the RHS of~\eqref{eq:amalgfpvN} and using
Corollary~\ref{thm:maincor}, we have $\pi_\tau(\Afr)''\cong L(\Fb_t)$, where $t$ is as in~\eqref{eq:tval}.
Compressing proves part~(5) of the theorem.
\end{proof} 

Since $A_s$ is a simple C$^*$--algebra, the GNS representation $\pi_{\tau\restrict_{A_s}}$ gives an embedding
of the Popa algebra $A_s$ as a weakly dense C$^*$--subalgebra of the interpolated free group factor $L(\Fb_s)$.
We have therefore proved Property~(I) from the introduction.

\section{Approximately Finite Dimensional Traces} 

In this section we prove another technical result which will be 
needed for the proof of the approximation properties
in the next section.  Namely, we will show that 
the canonical free product traces $\tau_{n}$ on the algebras $A(n)$ 
(notation as in the proof of Theorem \ref{thm:paconstruction}) are 
all weakly approximately finite dimensional in the sense of 
\cite[Definition 3.1]{brown:AFDtraces} (see Corollary 
\ref{thm:AFDtraces} below).  For the proof, we will need a lifting result,
whose proof uses some preliminary lemmas.

Let $\Mcal$ be a II$_1$--factor with tracial state $\tau$.
Let $\omega$ be a free ultrafilter on $\Nats$ and let
$\pi_\omega:\ell^\infty(\Nats,\Mcal)\to\Mcal^\omega$ be the quotient 
map onto the ultraproduct of $\Mcal$.
Denote by $\tau_\omega$ the trace on $\Mcal^\omega$ gotten by taking $\tau$ at the limit as $n\to\omega$,
and let $\sigma_n:\ell^\infty(\Nats,\Mcal)\to\Mcal$ be the evaluation map at the $n$th position:
$\sigma_n((a_k)_{k=1}^\infty)=a_n$.

\begin{lem}\label{lem:projclose}
Let $a=a^*\in\Mcal$ and let $p=E_a([1/2,\infty))$ be the spectral projection of $a$ for the set $[1/2,\infty)$.
Then
\begin{equation}\label{eq:pa}
\|p-a\|_2\le2\|a^2-a\|_2\;.
\end{equation}
\end{lem}
\begin{proof}
Let $\mu$ be $\tau$ of spectral measure of $a$.
Then
\begin{align*}
\|p-a\|_2^2&=\int_\Reals(1_{[1/2,\infty)}(t)-t)^2\dif\mu(t) \\
\|a^2-a\|_2^2&=\int_\Reals t^2(1-t)^2\dif\mu(t)\;.
\end{align*}
If $t<1/2$ then $t^2<4t^2(1-t)^2$, while if $t\ge1/2$ then $(1-t)^2\le4t^2(1-t)^2$;
the inequality~\eqref{eq:pa} follows.
\end{proof}

\begin{lem}\label{lem:projlift}
Let $p,q\in\Mcal^\omega$ be projections, with $p\le q$.
Suppose $Q=(Q_n)_{n=1}^\infty\in\ell^\infty(\Nats,\Mcal)$
is a projection such that $\pi_\omega(Q)=q$ and $\tau(Q_n)=\tau_\omega(q)$ for all $n\in\Nats$.
Then there is a projection $P=(P_n)_{n=1}^\infty\in\ell^\infty(\Nats,\Mcal)$ such that $P\le Q$,
$\pi_\omega(P)=p$ and $\tau(P_n)=\tau_\omega(p)$ for all $n\in\Nats$.
\end{lem}
\begin{proof}
Let $A=(A_n)_{n=1}^\infty\in\ell^\infty(\Nats,\Mcal)$ be such that $A^*=A$ and $\pi_\omega(A)=p$.
Replacing $A$ with $QAQ$, we may assume $A_n=Q_nA_nQ_n$ for all $n\in\Nats$.
Let $P_n'=E_{A_n}([1/2,\infty))$ be the spectral projection.
Then $P'=(P_n')_{n=1}^\infty$ is a projection in $\ell^\infty(\Nats,\Mcal)$, $P'\le Q$ and using Lemma~\ref{lem:projclose}
we get $\pi_\omega(P')=p$.
Hence $\lim_{n\to\omega}\tau(P_n')=\tau_\omega(p)$.
For every $n\in\Nats$, let $P_n\in\Mcal$ be a projection such that $P_n\le Q_n$, $\tau(P_n)=\tau_\omega(p)$ and
either $P_n\le P_n'$ or $P_n'\le P_n$.
Then $\lim_{n\to\omega}\|P_n-P_n'\|_2=0$ and the projection $P=(P_n)_{n=1}^\infty$ is as required.
\end{proof}

\begin{lem}\label{lem:pilift}
Let $v\in\Mcal^\omega$ be a partial isometry and suppose there are projections $P=(P_n)_{n=1}^\infty$
and $Q=(Q_n)_{n=1}^\infty$ in $\ell^\infty(\Nats,\Mcal)$ such that $\pi_\omega(P)=v^*v$, $\pi_\omega(Q)=vv^*$
and $\tau(P_n)=\tau(Q_n)$ for every $n\in\Nats$.
Then there is a partial isometry $V=(V_n)_{n=1}^\infty\in\ell^\infty(\Nats,\Mcal)$ such that
$\pi_\omega(V)=v$, $V^*V=P$ and $VV^*=Q$.
\end{lem}
\begin{proof}
Let $B=(B_n)_{n=1}^\infty\in\ell^\infty(\Nats,\Mcal)$ be such that $\pi_\omega(B)=v$.
Replacing $B$ by $QBP$, we may assume $B_n=Q_nB_nP_n$ for all $n\in\Nats$.
Let $B_n=U_n|B_n|$ be the polar decomposition.
Let $F_n=E_{|B_n|}([1/2,\infty))$ be the spectral projection.
Since $\pi_\omega(|B|)=v^*v$, using Lemma~\ref{lem:projclose} we get $\pi_\omega((F_n)_{n=1}^\infty)=v^*v$.
Thus $\pi_\omega((U_nF_n)_{n=1}^\infty)=v$.
Furthermore, since $F_n\le U_n^*U_n$, $U_nF_n$ is a partial isometry.
From $B_n=Q_nB_nP_n$ we obtain $U_n^*U_n\le P_n$ and $U_nU_n^*\le Q_n$.
Thus $F_n\le P_n$ and $U_nF_nU_n^*\le Q_n$.
Let $W_n\in\Mcal$ be a partial isometry such that $W_n^*W_n=P_n-F_n$ and $W_nW_n^*=Q_n-U_nF_nU_n^*$.
Let $V_n=U_nF_n+W_n$.
Then $V_n^*V_n=P_n$ and $V_nV_n^*=Q_n$.
Since $\pi_\omega((P_n-F_n)_{n=1}^\infty)=0$, we find $\lim_{n\to\omega}(\tau(P_n)-\tau(F_n))=0$ and therefore
$\lim_{n\to\omega}\|W_n\|_2=0$.
Thus $\pi_\omega((V_n)_{n=1}^\infty)=v$.
\end{proof}

\begin{prop}
\label{prop:lift}
Let $A$ be a finite dimensional C$^*$--algebra and suppose $\alpha:A\to\Mcal^\omega$ is a unital $*$--monomorphism.
Then there is a unital $*$--homomorphism $\alphat:A\to\ell^\infty(\Nats,\Mcal)$ such that $\alpha=\pi_\omega\circ\alphat$
and $\tau\circ\sigma_n\circ\alphat=\tau_\omega\circ\alpha$ for all $n\in\Nats$.
Furthermore, if $B\subseteq A$ is a C$^*$--subalgebra and if $\betat:B\to\ell^\infty(\Nats,\Mcal)$ is a $*$--homomorphism
such that $\pi_\omega\circ\betat=\alpha\restrict_B$ and $\tau\circ\sigma_n\circ\betat=\tau_\omega\circ\alpha\restrict_B$
for all $n\in\Nats$, then the $*$--homomorphism $\alphat$ above can be chosen so that
$\alphat\restrict_B=\betat$.
\end{prop}
\begin{proof}
We assume the existence of $B$ and the $*$--homomorphism $\betat$, since otherwise we can take $B=\Cpx1$.
Let $D\subseteq A$ be a maximal abelian subalgebra of $A$ such that $D\cap B$ is a maximal abelian subalgebra of $B$.
By successive application of Lemma~\ref{lem:projlift} to minimal projections of $D$, we can find a $*$--homomorphism
$\gammat:D\to\ell^\infty(\Nats,\Mcal)$, such that $\pi_\omega\circ\gammat=\alpha\restrict_D$,
$\gammat\restrict_{D\cap B}=\betat\restrict_{D\cap B}$ and
$\tau\circ\sigma_n\circ\gammat=\tau_\omega\circ\alpha\restrict_D$ for all $n\in\Nats$.

We will select some partial isometries $v_1,\ldots v_l\in A$ whose domain and range projections
are minimal projections in $D$
such that any choice of partial isometries $\vt_1,\ldots,\vt_l\in\ell^\infty(\Nats,\Mcal)$ satisfying
\begin{equation}\label{eq:vt}
\forall i\in\{1,\ldots,l\}\quad\pi_\omega(\vt_i)=v_i,\quad\vt_i^*\vt_i=\gammat(v_i^*v_i),
\quad\vt_i\vt_i^*=\gammat(v_iv_i^*)
\end{equation}
uniquely determines a $*$--homomorphism $\alphat:A\to\ell^\infty(\Nats,\Mcal)$ satisfying $\alphat(v_i)=\vt_i$,
$\alphat\restrict_B=\betat$ and $\alphat\restrict_D=\gammat$.
Once these partial isometries $v_1,\ldots,v_l$ have been found, the existence of $\vt_1,\ldots,\vt_l$ satisfying~\eqref{eq:vt}
is guaranteed by Lemma~\ref{lem:pilift}, and the $*$--homomorphism $\alphat$ will have been constructed.

Let $q_1,\ldots,q_m\in D\cap B$ be a maximal family of minimal projections in $B$ that are pairwise inequivalent in $B$.
For each $i\in\{1,\ldots,m\}$, let $p_i\in D$ be a minimal projection in $A$ such that $p_i\le q_i$.
Let $n(i)=\dim q_iD$ and let $w_{i,2},\ldots,w_{i,n(i)}$ be partial isometries in $A$ such that $w_{i,j}^*w_{i,j}=p_i$
and $q_i=p_i+\sum_{j=2}^{n(i)}w_{i,j}w_{i,j}^*$.
let $R_1,\ldots,R_{m'}$ be the equivalence classes of the set $\{p_1,\ldots,p_m\}$
under the relation of Murray--von Neumann equivalence.
Select a single element $p_{i_k}$ from each $R_k$, let $n'(k)=|R_k|$ and let $w_{k,2}',\ldots,w_{k,n'(k)}'$ be
partial isometries in $A$ such that $(w_{k,j}')^*w_{k,j}'=p_{i_k}$ and
$R_k=\{p_{i_k}\}\cup\{w_{k,j}'(w_{k,j}')^*\mid2\le j\le n'(k)\}$.
Now letting $v_1,\ldots,v_l$ be an enumeration of the set
\[
\{w_{i,j}\mid1\le i\le m,\,2\le j\le n(k)\}\cup\{w_{k,j}'\mid1\le k\le m',\,2\le j\le n'(k)\}\;,
\]
this collection of partial isometries has the desired property.
\end{proof}

\begin{cor}
\label{thm:AFDtraces} 
Let $A \cong M_{n} ({\Cpx})$ be a finite dimensional matrix algebra, 
$B \subset A$ be a unital subalgebra and $A*_{B}A$ denote the full 
amalgamated free product of $A$ with itself.  Then the canonical free 
product trace $\tau$ on $A*_{B}A$ is weakly approximately finite 
dimensional; i.e.\ there exists a sequence of unital, completely 
positive maps $\phi_{k} : A*_{B} A \to M_{m(k)}$ such that $\| 
\phi_{k}(ab) - \phi_{k}(a)\phi_{k}(b) \|_{2} \to 0$ and $\tr_{m(k)} 
\circ \phi_{k} (a) \to \tau (a)$, for all $a,b \in A*_{B}A$, where 
$\tr_{m(k)}$ denotes the normalized trace on $M_{m(k)}$ and $\| \cdot 
\|_{2}$ the induced 2-norm.    
\end{cor}

\begin{proof} The canonical free product trace 
$\tau$ on $A*_{B} A$ is defined as $\tau = \tr_{n} \circ\Et\circ 
\pi$ where $\pi : A*_{B}A \to (A,E) *_{B} (A,E)$ is the canonical 
quotient mapping onto the {\em reduced} amalgamated free product with 
respect to the $\tr_{n}$-preserving conditional expectation $E : A \to 
B$ and $\Et$ is the free product conditional expectation on the reduced amalgamated
free product C$^*$--algebra.

Since $(A,E) *_{B} (A,E)$ embeds into an interpolated free group 
factor (in a trace preserving way), by Theorem \ref{thm:RBR}, and 
since interpolated free group factors all embed into the ultrapower of 
the hyperfinite II$_{1}$ factor, $R^{\omega}$, it follows that we can 
find a unital $*$-homomorphism $\alpha : A*_{B}A \to R^{\omega}$ such 
that $\tau_{\omega}\circ\alpha = \tau$, where $\tau_{\omega}$ is the 
tracial state on $R^{\omega}$. 

It suffices to show that the $*$-homomorphism $\alpha$ 
lifts to a $*$-homomorphism $\beta : A*_{B}A \to l^{\infty}({\Nats}, R)$ (i.e.\ $\pi_{\omega} \circ \beta = \alpha$).  
Indeed, given such a $\beta$, one can compose the 
homomorphisms $\sigma_{n} \circ \beta : A*_{B}A \to R$ with 
conditional expectations onto larger and larger matrix subalgebras of 
$R$ to get the required maps $\phi_{k}$.  (Recall that $\sigma_{n} : 
l^{\infty}({\Nats}, R) \to R$ is defined by $\sigma_{n} ( (x_{i})_{i 
\in {\Nats}} ) = 
x_{n}$.) 

However, the existence of the desired $*$-homomorphism $\beta$ is
guaranteed by Proposition~\ref{prop:lift} and the proof is complete.
\end{proof}

\section{Approximation properties for the dense Popa algebras} 

\begin{thm}
\label{thm:mainthm}
For $1 < s < \infty$,
consider the interpolated free group factor $L(\Fb_s)$ with tracial state $\tau$
and let $L(\Fb_{s}) \subset B(H)$ be the corresponding GNS representation.
Then the finitely generated, weakly dense Popa 
algebra $A_{s} \subset 
L(\Fb_{s})$ constructed in~\S4 has the following properties:

\renewcommand{\labelenumi}{(\arabic{enumi})}
\begin{enumerate} 
    \item There exists a (non-normal) state $\varphi \in S(B(H))$ on 
    $B(H)$ such that $\varphi\restrict_{A_{s}} = \tau\restrict_{A_{s}}$ and $A_{s} 
    \subset B(H)_{\varphi} := \{ T \in B(H): \varphi(TS) = 
    \varphi(ST), \forall S \in B(H) \}$. 
    
    \item There exists a sequence of finite rank projections $P_{1}, 
    P_{2}, \ldots$ such that  
    \begin{enumerate}
	\item[(a)] $$\frac{\| [P_{n}, a] \|_{HS}}{\| P_{n} \|_{HS}} \to 0,$$ 
	\item[(b)] $$\frac{<aP_{n},P_{n}>_{HS}}{<P_{n}, P_{n}>_{HS}} 
	\to \tau(a),$$ 
    \end{enumerate}
     for all $a \in A_{s}$, where $<\cdot,\cdot>_{HS}$ (resp.\ $\| 
     \cdot \|_{HS}$) denotes the 
     Hilbert-Schmidt inner product (resp.\ norm) on finite rank operators.
    
    \item There exists a sequence of unital, completely positive 
    maps $\phi_{n} : A_{s} \to M_{k(n)}$ such that $\| 
    \phi_{n}(ab) - \phi_{n}(a)\phi_{n}(b) \|_{2} \to 0$ and $\tr_{k(n)} 
    \circ \phi_{n}(a) \to \tau(a)$, for all $a,b \in A_{s}$, where 
    $\tr_{k(n)}$ denotes the unique normalized trace on the $k(n)\times 
    k(n)$-matrices and $\| \cdot \|_{2}$ denotes the induced 2-norm. 
    
    \item There exists an idempotent, unital, completely postive 
    map $\Phi : B(H) \to L(\Fb_{s})$ such that $\Phi(a) = a$ 
    for all $a \in A_{s}$.  (Hence $\Phi(B(H))$ is a weakly dense, 
    injective operator sub-system of $L(\Fb_{s})$.) 
\end{enumerate}    
\end{thm}

\begin{proof}
Kirchberg (building on 
the celebrated work of Connes \cite{connes:classification}) has shown 
that the four properties are equivalent (cf.\ 
\cite[Theorem 3.6]{brown:AFDtraces}).
We will show part (4). 

It will again be convenient to identify each of the algebras $A(n)$ 
with their images in the algebra $\mathfrak{A}$ used in the proof of 
Theorem \ref{thm:paconstruction}.  Recall also that $A_{s} = 
p_{1}\mathfrak{A}p_{1}$ for a projection $p_{1} 
\in A(1) \subset \mathfrak{A}$ and that there exist projections 
$Q^{(n+1)} \in L(\Fb_{s}) \cap A(n)^{\prime}$ such that 
the weak closure of $Q^{(n+1)}A(n)$ is naturally isomorphic to 
$\pi_{\tau_{n}} (A(n))^{\prime\prime}$ (see part (p) in the proof of 
Theorem \ref{thm:paconstruction}).  

Since Corollary \ref{thm:AFDtraces} tells us that $\tau_{n}$ is a 
weakly approximately finite dimensional trace on $A(n)$ it follows 
that we can find a completely positive map $\Phi_{n} : B(H) \to 
Q^{(n+1)} A(n)^{\prime\prime}$, where $A(n)^{\prime\prime}$ denotes 
the weak closure of $A(n)$ in $L(\Fb_{s})$, such that 
$\Phi_{n} (x) = Q^{(n+1)}x$ for all $x \in A(n)$ (see \cite[Theorem 
3.6]{brown:AFDtraces}).   Taking any cluster point, in the topology of 
pointwise weak convergence, of the maps $\{ \Phi_{n} \}$ we get a 
unital, completely positive map $\Phi : B(H) \to L(\Fb_{s})$ 
such that $\Phi (a) = a$ for all $a \in A_{s}$.  One then replaces 
$\Phi$ with an idempotent such map by \cite[Theorem 
2.1]{blackadar:WEP} and the proof is complete. 
\end{proof}

We conclude this paper with a few remarks regarding Theorem 
\ref{thm:mainthm}.  First of all,
none of the properties (1)--(4) hold if $A_s$ is replaced by $L(\Fb_s)$, because
the interpolated free group factors are not hyperfinite.
Indeed, in this setting, each of the properties 
(1) - (4) in Theorem \ref{thm:mainthm} actually {\em characterizes} 
the hyperfinite II$_{1}$ factor -- i.e.\ $R$ is the unique II$_{1}$ 
factor which can be placed inside the centralizer of a state on 
$B(H)$ (property (1)) and is the unique II$_{1}$ factor satisfying 
Connes' F{\o}lner type condition (property (2)) and so on.

Secondly, as remarked in the introduction, property (4) shows that free group 
factors have the weak expectation property
of Lance~\cite{lance:nuclear} relative to $A_s$.

Finally, properties (1) - (4) are
almost never enjoyed by any sort of reduced free product 
C$^*$-algebra (taken with its GNS representation).
More precisely, no C$^{*}$-algebra $B$ which 
contains a unital copy of the reduced group C$^{*}$-algebra 
$C^{*}_{r} (\Fb_{2})$ has a weakly approximately finite dimensional tracial state,
which is the property described in~(3).
This is the case since existence of such a trace
clearly passes to subalgebras and it is  known 
that $C^{*}_{r} (\Fb_{2})$ 
has none (see \cite[Example 3.13]{brown:AFDtraces}).

\bibliographystyle{amsplain}

\begin{thebibliography}{99}

\bibitem{blackadar:WEP} B.~Blackadar, \emph{Weak expectations and 
injectivity in operator algebras}, Proc. Amer. Math. Soc. 
\textbf{68} (1978), 49-53. 
    
    
\bibitem{brown:AFDtraces} N.P. Brown, \emph{Tracial invariants, 
classification and II$_{1}$ factor representations of Popa algebras}, 
preprint 2001. 

\bibitem{brown:QDsurvey} 
\bysame, \emph{On quasidiagonal
C$^*$-algebras}, Proceedings of 1999 US-Japan conference on Operator
Algebras (to appear).

\bibitem{connes:classification} A.~Connes, \emph{Classification of
injective factors: cases II$_1$, II$_\infty$, III$_{\lambda}$,
$\lambda \neq 1$}, Ann. of Math. \textbf{104} (1976), 73--115.

\bibitem{D:fdim} K.~Dykema,
\emph{Free products of hyperfinite von Neumann algebras and free dimension,}
Duke Math. J. \textbf{69} (1993), 97--119.

\bibitem{D:interp} \bysame,
\emph{Interpolated free group factors,}
Pacific J.\ Math.\ {\bf 163} (1994), 123--135.

\bibitem{D:amalg} \bysame,
{\em Amalgamated free products of multi--matrix algebras and a construction of subfactors of a free group factor,}
Amer.\ J.\ Math.\ {\bf 117} (1995), 433-450.

\bibitem{D:subf} \bysame,
{\em Subfactors of free products of rescalings of a II$_1$--factor},
preprint (2002).

\bibitem{EL} R.~Exel, T.~Loring,
\emph{Finite-dimensional representations of free product $C\sp *$-algebras},
Internat. J. Math. \textbf{3} (1992), 469--476.

\bibitem{jung:hyperfinite} K.~Jung, \emph{The free entropy dimension 
of hyperfinite von Neumann algebras.}, preprint 2001.

\bibitem{lance:nuclear} E.C. Lance, \emph{On nuclear C$^*$-algebras},
J. Funct. Anal. \textbf{12} (1973), 157--176.

\bibitem{paulsen:cbmaps} V.~Paulsen, \emph{Completely bounded maps and
dilations}, Pitman Research Notes in Mathematics, vol. 146, Longman,
1986.


\bibitem{popa:simpleQD} S.~Popa, \emph{On local finite-dimensional
approximation of C$^*$-algebras}, Pacific J. Math. \textbf{181}
(1997), 141--158.


\bibitem{Ra94} F.\ R\u adulescu,
\emph{Random matrices, amalgamated free products and subfactors of
the von Neumann algebra of a free group, of noninteger index},
Invent. Math. \textbf{115} (1994), 347-389.

\bibitem{dvv:QDsurvey} D.V. Voiculescu, \emph{Around quasidiagonal
operators}, Integral Equations Operator Theory \textbf{17} (1993), 137--149.

\bibitem{V94} \bysame,
\emph{The analogues of entropy and of Fisher's information measure in free probability theory, II,}
Invent. Math. \textbf{118} (1994), 411--440.

\bibitem{dvv:partIII} \bysame, \emph{The analogues of entropy 
and of Fisher's information measure in free probability III: The 
absence of Cartan subalgebras}, Geom.
Funct. Anal. \textbf{6} (1996), 172--199.

\end{thebibliography}

\providecommand{\bysame}{\leavevmode\hbox to3em{\hrulefill}\thinspace}

\end{document}